\documentclass[11pt]{article}
\newcommand{\dated}{\mbox{} \hfill {\small [{\tt \today}]}} \usepackage{amsmath,amssymb,amscd}
\newcommand{\pf}[1]{\trivlist \item[\hskip\labelsep\it #1\ ]}
\newcommand{\varpf}[1]{\trivlist \item[\hskip\labelsep\sc #1:]}
\newcommand{\qedbox}{$\rlap{$\sqcap$}\sqcup$}
\newcommand{\qed}{\qquad \qedbox \endtrivlist}
\newcommand{\varqed}{\hfill \rule{0.6em}{0.6em} \endtrivlist}

\newenvironment{items}{
  \begin{enumerate} 
                    
  }{\end{enumerate}}

\newenvironment{keywords}{\noindent\small {\it Keywords\/}:}{\vskip 4pt}
\newenvironment{classification}{\noindent\small 2000 {\it Mathematics Subject
Classification\/}:}{\vskip 12pt}

\newcommand{\comps}{{\mathbb C}}

\newcommand{\tensor}{\otimes}

\newcommand{\A}{{\mathfrak A}}

\newtheorem{defi}{Definition}

\newenvironment{definition}{\begin{defi}\rm}{\end{defi}}
\newtheorem{theorem}{Theorem}

\newtheorem{lemma}{Lemma}

\newcommand{\SC}{\operatorname{\cal SC}_0}
\begin{document}
\title{Connes-amenability and normal, virtual diagonals \\ for measure algebras, II}
\author{{\it Volker Runde}\thanks{Research supported by NSERC under grant no.\ 227043-00.}}
\date{}
\maketitle
\begin{abstract}
We prove that the following are equivalent for a locally compact group $G$: (i) $G$ is amenable; (ii) $M(G)$ is Connes-amenable; (iii) $M(G)$ has a normal, virtual diagonal.
\end{abstract}
\section*{The result}
A Banach algebra $\A$ is called {\it dual\/} if it is a dual space such that multiplication in $\A$ is separately $w^\ast$-continuous. A Banach bimodule over a dual Banach algebra
is called {\it normal\/} if it is a dual space such that the module operations are separately $w^\ast$-continuous (see \cite{Run,LoA,Run2}).
\begin{definition}
A dual Banach algebra $\A$ is called {\it Connes-amenable\/} if, for every normal Banach $\A$-bimodule $E$, every $w^\ast$-continuous derivation $D \!: \A \to E$ is inner.
\end{definition}
\par
The notion of Connes-amenability was introduced for von Neumann algebras in \cite{JKR} (the name ``Connes-amenability'' seems to originate from \cite{Hel}); it is equivalent to a number
of important von Neumann algebraic properties such as injectivity, semidiscreteness, and being approximately finite-dimensional (see \cite[Chapter 6]{LoA} for a self-contained exposition and references
to the original literature).
\par
For arbitrary dual Banach algebras, Connes-amenability was first considered in \cite{Run}, and in \cite{Run2} it was shown that $M(G)$, the measure algebra of a locally compact group $G$, is
Connes-amenable if and only if $G$ is amenable.
\par
For any dual Banach algebra $\A$, let ${\cal L}^2_{w^\ast}(\A,\comps)$ denote the separately $w^\ast$-continuous bilinear functionals on $\A$. It is easy to see that the multiplication map
$\Delta \!: \A \tensor \A \to \A$ extends to ${\cal L}^2_{w^\ast}(\A,\comps)^\ast$ as a continuous $\A$-bimodule homomorphism.
\begin{definition}
Let $\A$ be a dual Banach algebra. A {\it normal, virtual diagonal\/} for $\A$ is an element ${\rm M} \in {\cal L}^2_{w^\ast}(\A,\comps)^\ast$ such that
\[
  a \cdot {\rm M} = {\rm M} \cdot a \quad\text{and}\quad a \Delta{\rm M} = a \qquad (a \in \A).
\]
\end{definition}
\par
If $\A$ has a normal, virtual diagonal, then it is Connes-amenable (\cite{CG,Eff}). The converse holds if $\A$ is a von Neumann algebra (\cite{Eff}). It is an open question --- likely with a negative answer ---
if Connes-amenability and the existence of normal, virtual diagonals are equivalent for arbitrary dual Banach algebras.
\par
In \cite{Run} and \cite{Run2}, we gave partial positive answers for $\A = M(G)$:
\begin{itemize}
\item If $G$ is compact, then $M(G)$ has a normal, virtual diagonal (\cite[Proposition 5.2]{Run} and \cite[Proposition 3.3]{Run2}).
\item If $G$ is discrete, then $M(G) = \ell^1(G)$ is Connes-amenable if and only if it has a normal, virtual diagonal (\cite[Corollary 5.4]{Run2}).
\end{itemize}
\par
In this note, we will prove (and thus extend \cite[Theorem 5.3]{Run2}):
\begin{theorem}
The following are equivalent for a locally compact group $G$:
\begin{items}
\item $G$ is amenable.
\item $M(G)$ is Connes-amenable.
\item $M(G)$ has a normal, virtual diagonal.
\end{items}
\end{theorem}
\section*{The proof}
For convenience, we quote the following well-known characterization of amenable locally compact groups (see \cite[Lemma 7.1.1]{LoA}, for example):
\begin{lemma}
A locally compact group $G$ is amenable if and only if there is a net $( f_\alpha )_\alpha$ of non-negative functions in the unit sphere of $L^1(G)$ such that
\begin{equation} \label{invlim} \tag{\mbox{$\ast$}}
  \sup_{x \in K} \| \delta_x \ast f_\alpha - f_\alpha \| \to 0
\end{equation}
for each compact subset $K$ of $G$.
\end{lemma}
\pf{Proof of the Theorem}
In view of \cite[Theorem 5.3]{Run2}, it is sufficient to show that (i) implies (iii).
\par
Let $( f_\alpha )_\alpha$ be a net as specified in the lemma. Define a net $( m_\alpha )_\alpha$ in $M(G \times G)$ by letting
\[
  \langle f, m_\alpha \rangle := \int_G f(x,x^{-1}) f_\alpha(x) \, dx \qquad (f \in {\cal C}_0(G \times G)),
\]
where $dx$ denotes integration with respect to left Haar measure on $G$. Let $y \in G$, and note that, for $f \in {\cal C}_0(G \times G)$,
\begin{eqnarray*}
  \langle f, (\delta_y \tensor \delta_e) \ast m_\alpha \rangle & = & \int_G f(yx,x^{-1}) f_\alpha(x) \, dx \\ 
    & = & \int_G f(x, x^{-1}y) f_\alpha(y^{-1}x) \, dx \qquad\text{(substitute $y^{-1}x$ for $x$)} \\
  & = & \int_G f(x,x^{-1}y)(\delta_y \ast f_\alpha)(x) \, dx
\end{eqnarray*}
and
\[
  \langle f, m_\alpha \ast (\delta_e \tensor \delta_y)\rangle = \int_G f(x,x^{-1}y) f_\alpha(x) \, dx.
\]
It follows from (\ref{invlim}) that
\begin{equation} \label{diaglim} \tag{\mbox{$\ast\ast$}}
  \sup_{y \in K} \| (\delta_y \tensor \delta_e) \ast m_\alpha -  m_\alpha \ast (\delta_e \tensor \delta_y) \| \to 0
\end{equation}
for each compact subset $K$ of $G$.
\par
By \cite[Proposition 3.1]{Run2}, we may identify the Banach $M(G)$-bimodules ${\cal L}^2_{w^\ast}(M(G),\comps)$ and 
\[
  \SC(G \times G) : = \{ f \in \ell^\infty(G \times G) : \text{$f( \cdot, x), f(x, \cdot) \in {\cal C}_0(G)$ for each $x \in G$} \}.
\]
Let $\cal U$ be an ultrafilter on the index set of $( m_\alpha )_\alpha$ that dominates the order filter. Define ${\rm M} \in \SC(G \times G)^\ast$ by letting
\[
  \langle f, {\rm M} \rangle := \lim_{\cal U} \int_{G \times G} f(x,y) \, dm_\alpha(x,y) \qquad (f \in \SC(G \times G))
\]
(since all functions in $\SC(G \times G)$ are measurable with respect to any Borel measure by \cite{Joh0}, the integrals do exist). It is routinely seen that $\Delta {\rm M} = \delta_e$.
\par
Let $\mu \in M(G)$ and let $f \in \SC(G \times G)$. Then we have:
\begin{eqnarray*}
  \lefteqn{| \langle f, \mu \cdot {\rm M} - {\rm M} \cdot \mu \rangle |} & & \\
  & = &  | \langle f \cdot \mu - \mu \cdot f , {\rm M} \rangle | \\
  & = & \left| \lim_{\cal U}  \int_{G \times G} \left( \int_G (f(zx,y) - f(x,yz)) \, d\mu(z) \right) dm_\alpha(x,y) \right| \\
  & = & \left| \lim_{\cal U}  \int_G \left( \int_{G \times G} (f(zx,y) - f(x,yz)) \, dm_\alpha(x,y) \right) d\mu(z) \right| \qquad\text{(by Fubini's theorem)} \\
  & \leq &  \lim_{\cal U}  \int_G \left| \int_{G \times G} (f(zx,y) - f(x,yz)) \, dm_\alpha(x,y) \right| d|\mu|(z) \\
  &  = & \lim_{\cal U}  \int_G \left| \int_{G \times G}  f(x,y) \, d((\delta_z \tensor \delta_e)\ast m_\alpha - m_\alpha \ast (\delta_e \tensor \delta_z))(x,y) \right|d|\mu|(z) \\
  & \leq & \lim_{\cal U} \int_G \| f \| \| (\delta_z \tensor \delta_e) \ast m_\alpha - m_\alpha \ast (\delta_e \tensor \delta_z) \| \, d|\mu|(z) \\
  & \to & 0 \qquad\text{(by (\ref{diaglim}) and the inner regularity of $|\mu|$)}.
\end{eqnarray*}
It follows that $\rm M$ is a normal, virtual diagonal for $M(G)$.
\qed
\dated
\vfill
\begin{tabbing}
{\it Address\/}: \= Department of Mathematical and Statistical Sciences \\
\> University of Alberta \\
\> Edmonton, Alberta \\
\> Canada T6G 2G1 \\[\medskipamount]
{\it E-mail\/}: \> {\tt vrunde@ualberta.ca} \\[\medskipamount]
{\it URL\/}: \> {\tt http://www.math.ualberta.ca/$^\sim$runde/runde.html}   
\end{tabbing} 
\end{document}